\documentclass[12pt]{article} 
 
\font\bigestssbx=cmssbx10 scaled \magstep4
\def\Asem#1#2{\mathop{\lower 5.5pt \hbox{\bigestssbx A}}_{#1}^{#2}}
\font\smallestssbx=cmssbx10 scaled \magstep2
\def\Asemsmall#1#2{\mathop{\lower 2.5pt \hbox{\smallestssbx A}}_{#1}^{#2}}
\def\asem#1#2{
              \ifmmode
                 \ifinner
                    \raise0.9pt\hbox{$\scriptstyle\bAs$}_{#1}^{#2}
                 \else
                    \Asem{#1}{#2}
                 \fi
              \fi
              }

\usepackage{amssymb}
\usepackage{amsfonts}
\usepackage{latexsym}
\usepackage{amsmath}
\usepackage{mathrsfs}

\usepackage{graphicx}
\usepackage{rotating}
\usepackage{psfrag}
\usepackage{wrapfig}
\usepackage{caption2}
\usepackage{epsfig}

\usepackage{color}

\usepackage{setspace}

\usepackage{xspace}
\usepackage{siunitx}
\usepackage[letterpaper]{geometry} 
\usepackage[numbers]{natbib}
\begin{document}


\title{A Monte Carlo packing algorithm for poly-ellipsoids and its comparison with packing generation using Discrete Element Model}

\author{
\small Boning Zhang, Eric B. Herbold, Richard A. Regueiro \\
\small Department of Civil, Environmental, and Architectural Engineering\\ 
\small University of Colorado, Boulder\\
\small Boulder, CO  80309
}


\maketitle

\begin{abstract}
Granular material is showing very often in geotechnical engineering, petroleum engineering, material science and physics. The packings of the granular material play a very important role in their mechanical behaviors, such as stress-strain response, stability, permeability and so on. Although packing is such an important research topic that its generation has been attracted lots of attentions for a long time in theoretical, experimental, and numerical aspects, packing of granular material is still a difficult and active research topic, especially the generation of random packing of non-spherical particles. To this end, we will generate packings of same particles with same shapes, numbers, and same size distribution using \textit{geometry} method and \textit{dynamic} method, separately. Specifically, we will extend one of Monte Carlo models for spheres to ellipsoids and poly-ellipsoids.

\vspace{.5cm}

\noindent
keywords: Discrete Element Method, Numerical Simulation, Numerical Method

\end{abstract}

\doublespacing

\clearpage

\section{Introduction}
Granular material is showing very often in geotechnical engineering, petroleum engineering, material science and physics \citep{buchalter1994, jaeger1996, Regueiro2014bio, Senseney2017, zhang2016thesis}. The packings of the granular material play a very important role in their mechanical behaviors, such as stress-strain response, stability, permeability and so on \citep{troadec1991, maggi2008, zhang2015, amirrahmat2018, Regueiro2014a, luo2022constitutive}. For example, \citet{troadec1991, 8329994, 10.1145/3543873.3587567} studied experimentally the effects of different packings on the stress-strain relation during compression of packed cylinders and found that the coordination number of the packing govern the stress-strain paths during compression. \citet{bassett2012,sadd2000,Jie2021BiddingVC} then further concluded experimentally and numerically, respectively, that force chain of packed granular materials plays significant effects on sound propagation in granular materials. And our previous researches \citep{Zhang2015fracture, la2016cumulative, ZHAO2018619} show that different packed sand grains behave differently in the simulations of Split Hopkinson Pressure Bar experiments.

Although packing is such an important research topic that its generation has been attracted lots of attentions for a long time in theoretical \citep{furnas1929,melissen1995,brouwers2006,torquato2009}, experimental \citep{furnas1929, bideau1993, 10.1145/3123266.3123296}, and numerical \citep{he1999, jia2001, jiang2003, dong2006, zhou2011, 7464796, peng2015probabilistic} aspects, packing of granular material is still a difficult and active research topic, especially the generation of random packing of non-spherical particles. There are two kinds of numerical generation methods of random packing based on \textit{geometry} and \textit{dynamic}, respectively. \textit{Geometry} method is moving and rotating particles in a packing by geometry constraints to reduce overlap between particles. After overlap is decreased to a small tolerance, the packing is viewed as a stable system and packing generation is finished. The \textit{geometry} method is efficient and has been applied successfully to generate random packings of monodisperse spheres, polydisperse spheres, ellipsoids, super-ellipsoids, and spherocylinders \citep{buchalter1994,yang1996,he1999,jia2001,kansal2002,fu2003,abreu2003,donev2007,maggi2008,stafford2010,delaney2010}. Among so many \textit{geometry} packing generation algorithms, Monte Carlo model is very famous and extended to many different variables. The \textit{dynamic} method is simulating the physical interactions between particles during packing process \citep{cheng2000,jiang2003,zhou2011,jiecheng-thesis, nia2014streaming} based on Discrete Element Method (DEM). DEM \citep{Cundall1979} is used to simulate the motion of individual particles/grains within an overall deformation and/or flow of a granular medium, which is borrowed to simulate the dynamic process during packing generation. For example, \citet{cheng2000,jiemulti} randomly generated circular particles inside a box without any overlaps by removing the overlapped particles, then these random particles are settled down by gravity force to generate random packings. \citet{jiang2003} generated particles randomly in a container without any overlaps, and then compress these particles to reach specific packing fraction. \citet{zhou2011} generated packings by pouring ellipsoidal particles from a certain height into a container, which is more physically related to the generation of packings in experiments. These descriptions of the two methods suggest that \textit{geometry} methods are much more faster than \textit{dynamic} methods, since it takes a lot of time for the particles in the latter to come to rest. While \textit{dynamic} methods considered the interactions between particles mechanically, and the pouring method is more similar physically to the generation of packings in experiments. Although both \textit{geometry} methods and \textit{dynamic} methods can generate packings with similar fractions as in experiments, there are worries about what it will influence on the packings if \textit{geometry} methods do not consider packing as a dynamic process involving forces between particles \citep{jia2001,dong2006,an2008effect,zhou2011}. The motivation of this paper is then to explore the influence and to answer if we can take advantages of the efficiency of \textit{geometry} methods without introducing much difference to the final packings of particles compared with \textit{dynamic} generating methods. To this end, we will generate packings of same particles with same shapes, numbers, and same size distribution using \textit{geometry} method and \textit{dynamic} method, separately. Specifically, we will extend one of Monte Carlo models proposed by \citet{he1999} for spheres to ellipsoids and poly-ellipsoids \citep{Peters2009}. To simulate the dynamic process during packings, we will adopt gravitational deposition \citep{Yan2010} of particle assembly which is similar as the pouring method described in \citep{zhou2011}. After packings are generated by the two methods separately, the corresponding packings will be analyzed and compared with respect to packing metrics, such as packing fraction, Coordination Number (CN), Radia Distribution Function (RDF) \citep{konakawa1990,he1999,zhou2011}, fabric tensor \citep{Yimsiri2010}, and also force chain \citep{sadd2000,bassett2012}. 

Currently, most of packings are generated for spheres and ellipsoids. Some researchers studied and concluded that particle shape plays a crucial role in packings of ellipsoids \citep{bezrukov2006}, super-ellipsoids \citep{delaney2010}, hyperspheres \citep{skoge2006packing}, and sphereocylinders \citep{abreu2003}. Thus this paper also considered packings of different shapes generated by the two methods, such as spheres, ellipsoids (prolate and oblate) \citep{bezrukov2006}, and poly-ellipsoids (carrot and half-dome) \citep{Peters2009}. With poly-ellipsoids, we can create many unsymmetrically shaped particles, it is very important since most of geotechnical material are not regular or symmetric, and no one else has studied packings of unsymmetric particles.

\section{Algorithms}
This paper considered two packing generation methods, an extended Monte Carlo Model and DEM simulation of gravitational deposition, which are belonging to the two major packing generation methods in community, respectively. The two methods will be used to generate packings of same particle assemblies, and then the packings will be analyzed and compared to study the effects caused by not considering dynamic interactions between particles in Monte Carlo packing algorithm.

\subsection{Extended Monte Carlo Model}
\citet{he1999} presented a variable of Monte Carlo model to generate random packings of unequal spherical particles obeying any specified distributions. \citet{he1999} first generated a certain number ($ n $) of spherical particles following a specified size distribution, and these particles are uniformly randomly placed within a cubic domain with the initial size as $ L_0 $, which is determined as \citep{he1999}
\begin{align}
L_0 = \left[\frac{1}{\Phi_0}\sum_{i=1}^n v_i\right]^{1/3} \label{eqn_L0}
\end{align}
where $ \Phi_0 $ is an arbitrary initial packing density which can be unphysically higher than random close packing fraction, in this paper $ \Phi_0 $ is taken as $ 0.86 $ following \citet{he1999}'s suggestion, and $ v_i $ is the volume of particle $ p_i $. The meaning of Eqn.\ref{eqn_L0} is generating a cubic domain containing the initial particle assembly with a packing fraction as $ \Phi_0 $. 

\begin{figure}[!h]
\begin{center}
\includegraphics[scale=0.6]{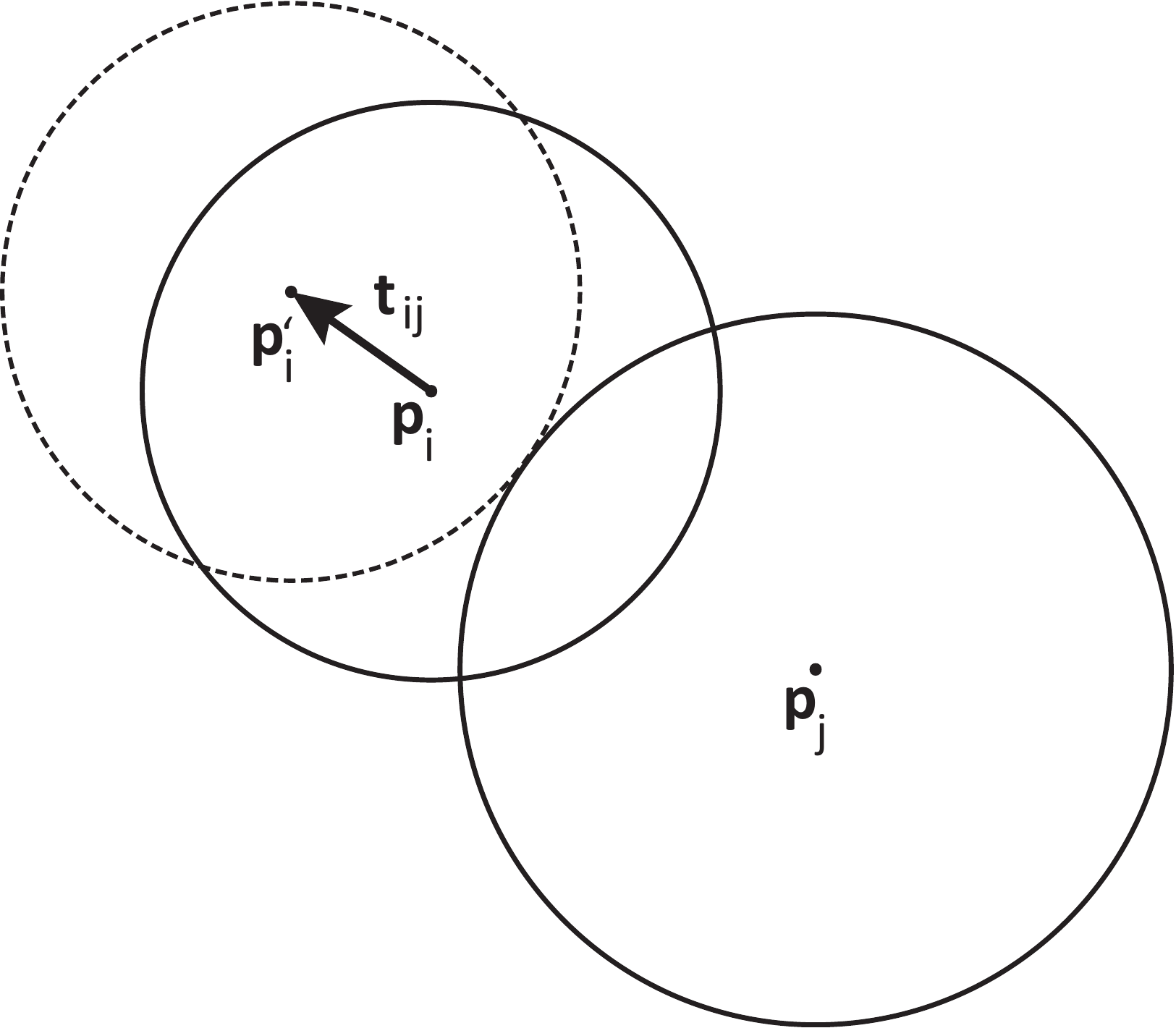}
\end{center}
\caption{\small Push outward particle $ p_i $ in the direction of branch vector connecting the centroids of particle $ p_i $ and $ p_j $ by translation vector $ \textbf{t}_{ij} $ in order to eliminate the overlap between particle $ p_i $ and $ p_j $, motivated by \citet{he1999}.}
\label{sphere_tij}
\end{figure}

There will be lots of overlaps between particles in this initial packing. During each step, \citet{he1999} moved an overlapped particle by the sum vector considering all overlaps of this particle. For example, the new position of particle $ p_i $ with $ n_i $ overlapped particles is given by \citep{he1999} 
\begin{align}
\textbf{p}_i^{\prime} = \textbf{p}_i+\frac{1}{n_i}\sum_{j=1}^{n_i}\textbf{t}_{ij}
\end{align}
where $ \textbf{t}_{ij} $ is the translation vector of particle $ p_i $ caused by the overlapping with particle $ p_j $, which is shown in Fig.\ref{sphere_tij}. After a certain number of iterations, the mean of the relative overlap will eventually drop to a very small preset tolerance, such that the packing is regarded as a stable packing.

\begin{figure}[!h]
\begin{center}
\includegraphics[scale=1.0]{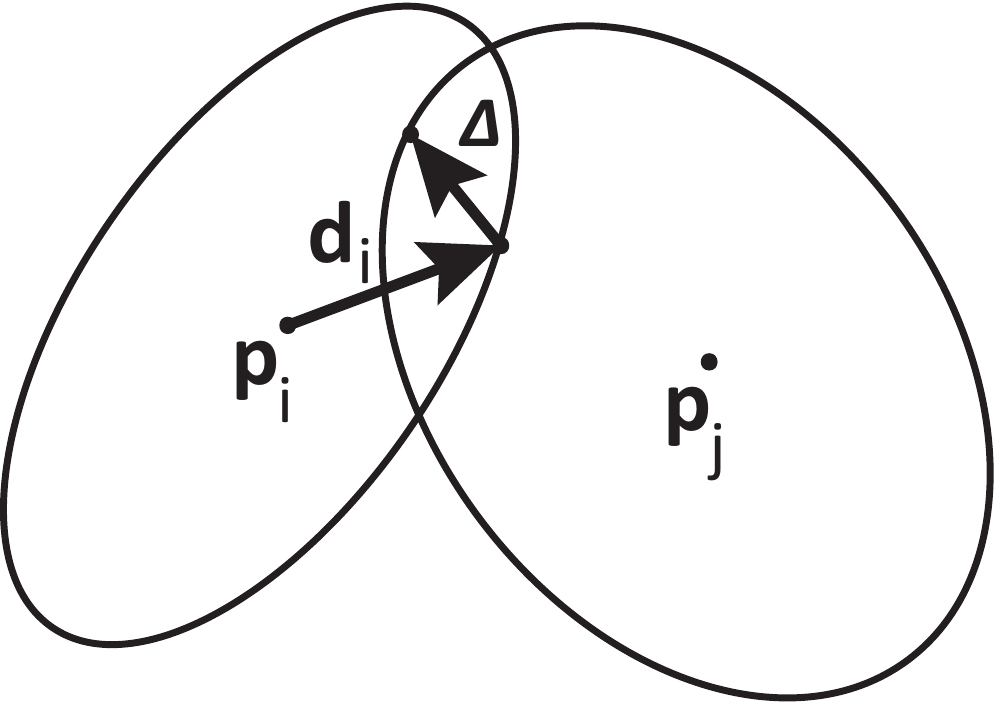}
\end{center}
\caption{\small $2$D illustration of overlapped particles. $ \boldsymbol{\Delta} $ is the overlap vector connecting two contact points $ i $ and $ j $, which are on the surface of particle $ p_i $ and $ p_j $. $ \textbf{d}_i $ is pointing from centroid $ p_i $ to point $ i $.}
\label{ellipsoid_trans_rotate}
\end{figure}

To extend this Monte Carlo model for ellipsoids and poly-ellipsoids, the major changes are the contact detection and rotations of overlapped particles. The contact detection algorithm is to find contacts between particles and their overlaps. For spheres, the contact detection algorithm is trivial, however it is not straightforward for ellipsoids and poly-ellipsoids, the contact detection algorithms for which are described in \citet{Yan2010} and \citet{zhang2016construction}, respectively. The other change is to add rotations of overlapped particles for ellipsoids and poly-ellipsoids. Suppose we have two overlapping particles ($ p_i $ and $ p_j $) as in Fig.\ref{ellipsoid_trans_rotate}, then the translation and rotation of particle $ p_i $ caused by particle $ p_j $ in our model are, respectively
\begin{align}
\textbf{t}_{ij} &= \frac{m_i}{m_i+m_j}\boldsymbol{\Delta} \\
\textbf{r}_{ij} &= \textbf{d}_i \times \textbf{t}_{ij}
\end{align}
where $ m_i $, $ m_j $ are the mass for particles $ p_i $ and $ p_j $, and $ \boldsymbol{\Delta} $, $ \textbf{d}_i $ are shown in Fig.\ref{ellipsoid_trans_rotate}. $ \boldsymbol{\Delta} $ is the overlap vector pointing from point $ i $ to point $ j $, which are contact points on the surface of particle $ p_i $ and $ p_j $, respectively. $ \textbf{d}_i $ is a vector pointing from centroid $ p_i $ to point $ i $. Similarly, if particle $ p_i $ has $ n_i $ overlapped particles, then the overall translation vector and rotation vector are 
\begin{align}
\textbf{T}_{ij} &= \sum_{j=1}^{n_i} t_{ij} \\
\textbf{R}_{ij} &= \sum_{j=1}^{n_i} r_{ij}
\end{align}
Then after a certain number of steps, the mean of relative overlap will drop below a tolerance and the packing generation is finished. In this research, $ 1\times 10^{-4} $ is accepted as the tolerance for the mean of relative overlap.

\subsection{DEM simulation of gravitational deposition}

DEM \citep{Cundall1979} is often used in geotechnical engineering to simulate the motion of individual particles/grains within a deforming/flowing granular medium. According to \citet{Cundall1979}, the governing equation for the translation and rotation of DEM particle $ i $ in an assembly are
\begin{equation}
\begin{aligned}
m_i \ddot{\textbf{u}}_i &= \textbf{F}_i \\
I_i \ddot{\boldsymbol{\theta}}_i &= \textbf{M}_i
\end{aligned} \label{eqn_govern_DEM}
\end{equation}
where $ \textbf{u} $ is the particle displacement; $ \boldsymbol{\theta} $, the orientational vector of the particle; $ m $, the particle mass; $ I $, the moment of inertia of the particle; \textbf{F}, the resultant force, which includes body force and contact forces by its overlapping particles; and $ \textbf{M} $, the resultant moment about the principal axes of inertial frame. In this DEM model, Hertz-Mindlin contact theory \citep{mindlin1949} is applied which includes nonlinear elasticity and slip, with addition of Coulomb friction model to evaluate tangential stick-slip conditions. Central difference time integration method is used to solve Eqn.\ref{eqn_govern_DEM}.

To generate packings of granular material, a gravitational deposition of particle assembly is simulated dynamically using DEM. The particles will be positioned inside a container with open top boundary. The particles are positioned initially without any overlaps with other particles and boundaries. Then these particles will be dropped by the gravity forces at zero initial velocities. After all the particles come to rest, the packing generation is completed.


\subsection{Initial particle size distribution}
\begin{figure}[!h]
\begin{center}
\includegraphics[scale=0.6]{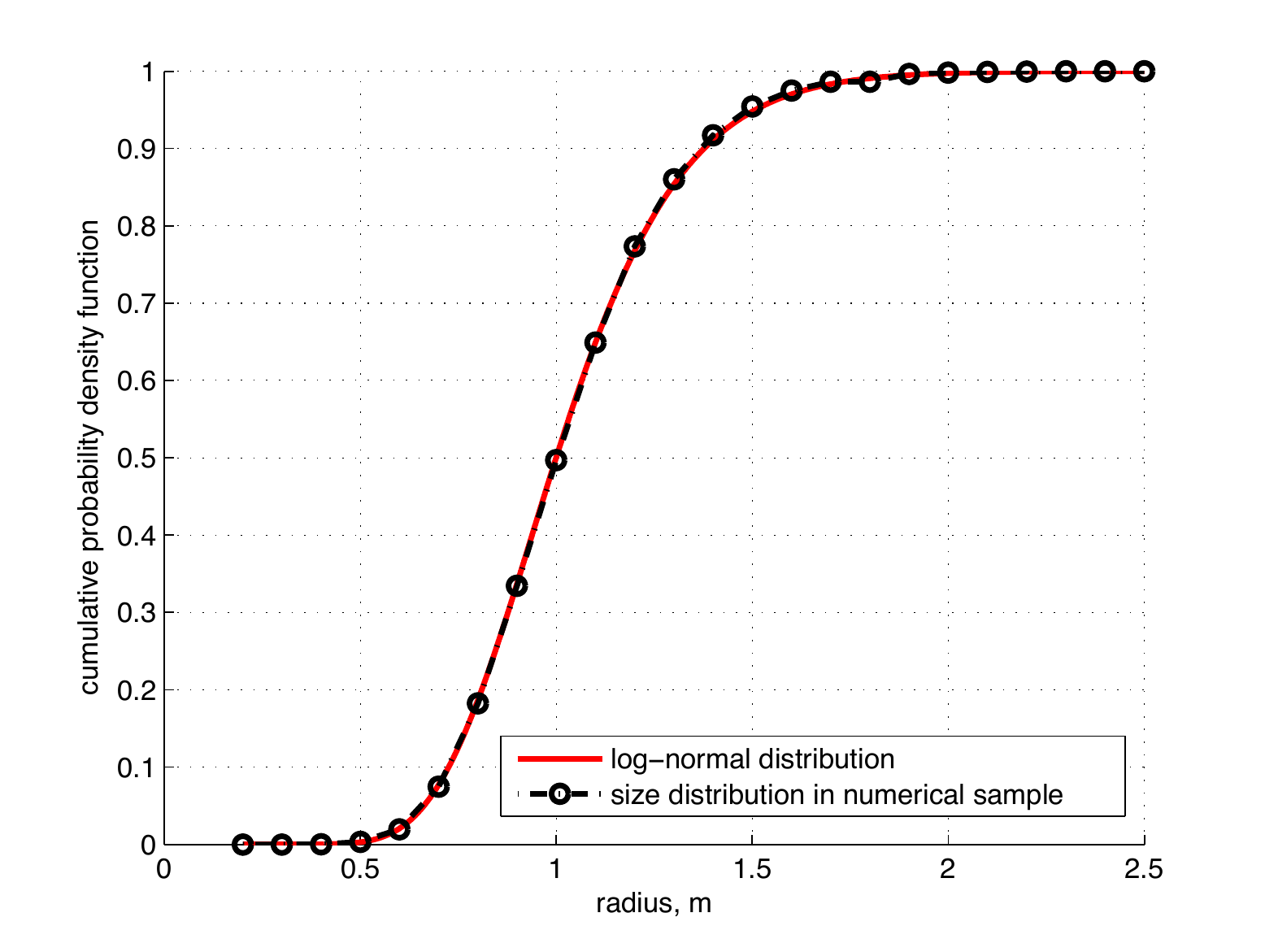}
\end{center}
\caption{\small Particle size distribution in numerical sample is compared well with the analytical log-normal distribution.}
\label{size_distribution}
\end{figure}

Similar as \citet{he1999}, sizes of particles obey truncated log-normal distribution, the probability density function of particle with radius $ r $ is 
\begin{align}
f(r) = \frac{1}{\sqrt{\left(2\pi\right)}\sigma r}e^{-\left(ln r-ln r_0\right)^2/\left(2\sigma^2\right)} \label{eqn_distr}
\end{align}
In this research, the mean radius is chosen as $ r_0=\SI{1}{\meter} $, and the standard deviation $ \sigma = \SI{0.25}{\meter} $, and  we only allow particle radius between $ \SI{0.2}{\meter} $ and $ \SI{2.5}{\meter} $. The particle size distribution generated in numerical sample is compared with analytical log-normal distribution as shown in Fig.\ref{size_distribution}.

There are five types of shaped particles in this research, i.e. sphere, prolate, oblate, carrot, and half-dome. Prolate and oblate are ellipsoids, while carrot and half-dome are poly-ellipsoids. For spheres, their radii are determined by Eqn.\ref{eqn_distr}. For ellipsoids and poly-ellipsoids, we make their major semi-lengths equal to $ r $ determined by Eqn.\ref{eqn_distr}. Particle assembly with each shaped particle will be only generated once and then the same assembly will be used in the extended Monte Carlo model and DEM gravitational deposition.




\newpage
\bibliographystyle{plainnat}
\bibliography{./references}






\end{document}